\documentclass[12pt]{amsart}
\usepackage{amsfonts,amssymb,amsmath}
\usepackage[english]{babel}
\usepackage{amsfonts,amssymb}
\usepackage{geometry}
\newtheorem{theorem}{Theorem}
\newtheorem{lemma}{Lemma}
\newtheorem*{baslem*}{Basic Lemma}

\theoremstyle{definition}
\newtheorem{definition}{Definition}
\theoremstyle{remark}
\newtheorem{remark}{Remark}

\theoremstyle{definition}
\newtheorem*{definition*}{Definition}
\theoremstyle{remark}
\newtheorem*{remark*}{Remark}
\newtheorem*{example*}{Example}

\begin{document}

\title{Regularization of singular Sturm-Liouville equations}

\author{Andrii Goriunov, Vladimir Mikhailets}

\address{Institute of Mathematics of National Academy of Sciences of Ukraine \\
         Tereshchenkivska str., 3 \\
         Kyiv-4 \\
         Ukraine \\
         01601}

\email[Andrii Goriunov]{goriunov@imath.kiev.ua}
\email[Vladimir Mikhailets]{mikhailets@imath.kiev.ua}

\subjclass[2010]{Primary 34L40; Secondary 34B08, 47A10}

\begin{abstract}
 The paper deals with the singular Sturm-Liouville expressions
  $$l(y) = -(py')' + qy$$
  with the coefficients
  $$q = Q', \quad 1/p,
  Q/p, Q^2/p \in L_1,
  $$
  where the derivative of the function $Q$ is
  understood in the sense of distributions. Due to a new
  regularization, the corresponding operators are correctly defined
  as quasi-differentials. Their resolvent approximation is
  investigated and all self-adjoint and maximal dissipative
  extensions and generalized resolvents are described in terms of
  homogeneous boundary conditions of the canonical form.
\end{abstract}

\keywords{Sturm-Liouville problem, quasi-differential expression, singular coefficients, resolvent approximation,
self-adjoint extension, generalized resolvent}

\maketitle
\section{Introduction}
This paper studies operators generated by the differential expressions
\begin{equation}
\label{vyraz} l(y) = -(py')'(t) + q(t)y(t), \quad t \in \mathcal{J}
\end{equation}
on a finite interval $\mathcal{J} := (a, b)$.

If the coefficients in (\ref{vyraz}) are real-valued and
\begin{equation}
q \in C(\overline{\mathcal{J}}),  \quad  0 < p \in C^1(\overline{\mathcal{J}}),
\end{equation}
then the equation $l(y) = f$ is a differential Sturm-Liouville
equation that has been investigated quite comprehensively. A modern
exposition of the classical Sturm-Liouville theory may be found in
many studies.  Principal statements of this theory remain true under
the weaker assumptions
\begin{equation}\label{Z cond}
q, 1/p \in L_1\left(\mathcal{J}, \mathbb{C} \right) = :L_1,
\end{equation}
see \cite{Z} and references therein.  This is achieved through
a regularization of the expression $l(y)$ applying Shin-Zettl
quasi-derivatives.  They were introduced in \cite{Shin} and later
generalized in \cite{Zettl}, see also \cite{EM}.

A further essential development of that approach was achieved in the
paper~\cite{S-Sh}.  It was proved there that if $p(t) \equiv 1$,
then the condition on $q$ may be significantly weakened.  Namely, it
is sufficient to suppose that
\begin{equation}\label{S-Sh cond}
p(t) \equiv 1, \quad q = Q', \quad Q \in L_2\left(\mathcal{J}, \mathbb{C}\right) =: L_2,
\end{equation}
where the derivative of the function $Q$ is understood in the sense
of distributions.  Note that the one-dimension Schr\"{o}dinger
operators with potentials that are Radon measures were introduced
and investigated long before that by physicists applying operator
theory methods (see \cite{Albeverio} and references therein).

The main goal of this paper is to define and investigate
Sturm-Liouville operators on a finite interval $\mathcal{J}$ under
the assumptions more general than those in~(\ref{Z cond}) and
(\ref{S-Sh cond}),
\begin{equation}\label{GM cond}
q = Q', \quad 1/p, Q/p, Q^2/p \in L_1.
\end{equation}

To achieve this goal, in Section 2, we propose a new regularization
of the formal differential expression (\ref{vyraz}) under
assumptions (\ref{GM cond}) by means of Shin-Zettl
quasi-derivatives.  We also define the corresponding maximal and
minimal operators on the Hilbert space $L_2$.  If conditions (\ref{Z
  cond}) hold, then these operators coincide with the classical ones
and, under assumptions (\ref{S-Sh cond}), they are identical to the
operators introduced in \cite{S-Sh}.

Section 3 shows that, in the case of two-point boundary conditions,
resolvents of the constructed operators may be approximated in the
sense of the norm with resolvents of other Sturm-Liouville
operators; for instance, ones that have more regular coefficients.

In Section 4, the minimal operator is supposed to be symmetric and all
its self-adjoint extensions are described in terms of the
homogeneous boundary conditions of the canonical form.

In addition, in Section 5, all maximal dissipative extensions and
generalized resolvents of the minimal symmetric operator are
described in the same form.

Extensions in Sections 4 and 5 are described by applying the boundary
triplet theory (see \cite{Gorbachuk} and references therein).  They
are parametrized by certain classes of operators on $\mathbb{C}^2$,
and this parametrization is bijective and continuous.  Also,
separated boundary conditions are singled out.

Note that in the case where $p(t) \equiv 1$, the results of Sections
3 and 4 improve the corresponding results of \cite{S-Sh} where
stronger conditions are required for the approximation, and self-adjoint
extensions are described on the basis of the Glazman-Krein-Naimark
theory.  Results of Section 5 deal with the questions not considered
in \cite{S-Sh}.
\section{Regularization of singular expression}
Consider the formal differential expression (\ref{vyraz}), assuming that
conditions (\ref{GM cond}) hold.  We introduce the quasi-derivatives
\begin{align*}
&D^{[0]} y = y, \\
&D^{[1]} y = py' - Qy, \\
&D^{[2]} y = (D^{[1]} y)' + {Q\over p}D^{[1]} y + {Q^2\over p}y.
\end{align*}
Then expression (\ref{vyraz}) is defined to be the quasi-differential expression
$$l[y] := -  D^{[2]} y.$$
\begin{definition}
  A solution of the Cauchy problem for the resolvent equation
\begin{equation}\label{cauchy pr 1}
l[y] - \lambda y = f\in L_2, \quad y(c) = \alpha_1, \quad
(D^{[1]}y)(c) = \alpha_2,
\end{equation}
where $c \in \overline{\mathcal{J}}$ and $\alpha_1, \alpha_2$ are
arbitrary complex numbers, is defined to be the first component of
the solution of the Cauchy problem for the correspondent system of
the first order differential equations
\begin{equation}\label{cauchy pr 2}
w'(t)=A_\lambda(t)w(t) + \varphi(t), \quad w(c) = (\alpha_1,
\alpha_2),
\end{equation}
where $w(t) = (y(t), D^{[1]}y(t))$,
the matrix-valued function is
$$
A_\lambda(t):=\left ( \begin{array}{cc}
\frac{Q}{p}& \frac{1}{p}\\
-\frac{Q^2}{p} - \lambda &-\frac{Q}{p}
\end{array}\right) \in L_1^{2\times 2},
$$
and $\varphi(t) := (0, -f(t))$.
\end{definition}

\begin{lemma}\label{cauchy pr unique}
  Problem (\ref{cauchy pr 1}), with assumptions (\ref{GM cond}), has
  only a unique solution defined on $\overline{\mathcal{J}}$.
\end{lemma}

\begin{proof}[Proof of Lemma \ref{cauchy pr unique}]
  Problem (\ref{cauchy pr 2}) with $A_\lambda(\cdot) \in
  L_1^{2\times 2}$ has only a unique solution for any $c \in
  \overline{\mathcal{J}}$ and $\left(\alpha_1, \alpha_2\right) \in
  \mathbb{C}^2$ due to Theorem 1.2.1 \cite{Z}.  This implies the
  statement of Lemma \ref{cauchy pr unique} by Definition 1.
\end{proof}

The quasi-differential expression $l[y]$ gives rise to the \emph{maximal}
quasi-differential operator
 $$ L_{\operatorname{max}}:y \to
l[y],\quad \operatorname{Dom}(L_{\operatorname{max}}) := \left\{y
  \in L_2 \left| y, D^{[1]}y \in AC(\overline{\mathcal{J}},
    \mathbb{C}), D^{[2]} y \in L_2\right.\right\}
$$
on the Hilbert space $L_2$ (see \cite{Zettl, EM}).  The
\emph{minimal} quasi-differential operator is defined as a
restriction of the operator $L_{\operatorname{max}}$ onto the set
$$ \operatorname{Dom}(L_{\operatorname{min}})  := \left\{y \in
\operatorname{Dom}(L_{\operatorname{max}}) \left| D^{[k]}y(a) =
D^{[k]}y(b) = 0, k = 0,1\right.\right\}.$$

\begin{remark}
  One can easily see that if $Q$ is replaced with $\widetilde{Q} :=
  Q + c$, $c \in \mathbb{C}$, then the operators
  $L_{\operatorname{max}}, L_{\operatorname{min}}$ do not change.

  If the coefficients in (\ref{vyraz}) satisfy (\ref{Z cond}), then the
  operators $L_{\operatorname{max}}, L_{\operatorname{min}}$
  introduced above coincide with the usual maximal and minimal
  Sturm-Liouville operators \cite{Z}.
\end{remark}

Consider the expression$$
l^+(y) =
-(\overline{p}y')'(t) + \overline{q}(t)y(t),
$$
formally adjoint to (\ref{vyraz}), where the bar denotes complex
conjugation. Denote by $L^+_{\operatorname{max}}$ and
$L^+_{\operatorname{min}}$ the maximal and the minimal operators
generated by this expression on the space $L_2$. Then  results of
this section, together with results of \cite{EM} for general
quasi-differential expressions, yield following theorem.

\begin{theorem} \label{L adjoint} The operators
  $L_{\operatorname{min}}$, $L^+_{\operatorname{min}}$,
  $L_{\operatorname{max}}$, $L^+_{\operatorname{max}}$ are closed
  and densely defined on the space $L_2$,
 $$L_{\operatorname{min}}^* = L^+_{\operatorname{max}},\quad L_{\operatorname{max}}^* = L^+_{\operatorname{min}}.$$
 In the case where $p$ and $Q$ are real-valued, the operator
 $L_{\operatorname{min}} = L^+_{\operatorname{min}}$ is symmetric
 with the deficiency index $\left({2,2} \right)$, and
$$L_{\operatorname{min}}^* = L_{\operatorname{max}},\quad
L_{\operatorname{max}}^* = L_{\operatorname{min}}.$$
 \end{theorem}

\section{Approximation of resolvent}
Consider the class of quasi-differential expressions $l_\varepsilon
[y] = -D_\varepsilon^{[2]}y$ with the coefficients
$$p_\varepsilon, q_\varepsilon = Q'_\varepsilon, \quad \varepsilon \in [0, \varepsilon_0].$$

On the Hilbert space $L_2$ with norm $\|\cdot\|_2$, these
expressions generate the operators
$L^\varepsilon_{\operatorname{min}} $,
$L^\varepsilon_{\operatorname{max}}$ for every $\varepsilon$. Let
$\alpha(\varepsilon),\beta(\varepsilon)\in \mathbb{C}^{2\times 2}$
be matrices and consider the vectors
$$\mathcal{Y}_\varepsilon(a):=\left\{y(a),D^{[1]}_\varepsilon
y(a)\right\},\quad
\mathcal{Y}_\varepsilon(b):=\left\{y(b),D^{[1]}_\varepsilon
y(b)\right\} \in \mathbb{C}^2.$$

Consider the quasi-differential operators
$$L_\varepsilon y =
l_\varepsilon[y],\quad \operatorname{Dom}(L_\varepsilon) =
\left\{\left.y \in
\operatorname{Dom}\left(L^\varepsilon_{\operatorname{max}}\right)\right|
\alpha(\varepsilon)\mathcal{Y}_\varepsilon(a)+\beta(\varepsilon)\mathcal{Y}_\varepsilon(b)=0\right\}.$$

It is evident that $L^\varepsilon_{\operatorname{min}} \subset
L_\varepsilon \subset L^\varepsilon_{\operatorname{max}}, \quad
\varepsilon \in [0, \varepsilon_0].$

We denote by $\rho(L)$ the resolvent set of the operator $L$.
Recall that the operators $L_\varepsilon$ converge to the operator
$L_0$ in the sense of the norm resolvent convergence, $L_\varepsilon
\stackrel{R}{\rightarrow} L_0$, if the there is a number $\mu \in
\mathbb{C}$ such that $\mu \in \rho(L_0)$ and $\mu \in
\rho(L_\varepsilon)$ (for all sufficiently small $\varepsilon$), and
$$\|(L_\varepsilon - \mu)^{-1} - (L_0 - \mu)^{-1}\| \rightarrow 0, \quad \varepsilon \rightarrow 0+.$$
This definition does not depend on the choice of the point $\mu \in \rho(L_0)$ \cite{K}.

In the case where the matrices
$\alpha(\varepsilon),\beta(\varepsilon)$ do not depend on
$\varepsilon$ and $p_\varepsilon(t) \equiv 1$, it is shown
in~\cite{S-Sh} that if $\|Q_\varepsilon - Q_0\|_2 \rightarrow 0$ for
$\varepsilon \rightarrow 0+$ and the resolvent set of the operator
$L_0$ is not empty, then $L_\varepsilon \stackrel{R}{\rightarrow}
L_0$. The following theorem generalizes this result.

\begin{theorem}\label{res conv part}
  Suppose $\rho(L_0)$ is not empty and, for $\varepsilon \rightarrow
  0+$, the following conditions hold:
\begin{align*}
(1)&\,\, \|1/p_\varepsilon - 1/p_0\|_1 \rightarrow 0,\\ (2)&\,\,
\|Q_\varepsilon/p_\varepsilon - Q_0/p_0\|_1 \rightarrow 0,\\
(3)&\,\, \|Q^2_\varepsilon/p_\varepsilon - Q^2_0/p_0\|_1
\rightarrow 0,\\ (4)&\,\,\alpha(\varepsilon)\rightarrow
\alpha(0),\quad \beta(\varepsilon)\rightarrow\beta(0),
\end{align*}
where $\|\cdot\|_1$ is the norm in the space $L_1(\mathcal{J},
\mathbb{C})$.  Then $L_\varepsilon \stackrel{R}{\rightarrow} L_0$.
\end{theorem}

\begin{remark}
  In the case where $p_\varepsilon(t) \equiv 1$, condition (1) is
  automatically fulfilled and conditions (2) and (3) are weaker than
  the assumption that $\|Q_\varepsilon - Q_0\|_2 \rightarrow 0$.
\end{remark}

To prove Theorem \ref{res conv part} we will need some auxiliary
results.

We start by introducing the following definition (\cite{MR1, MR2}).

\begin{definition}
  Denote by $\mathcal{M}^{m}(\mathcal{J}) =:\mathcal{M}^{m},$ $ m\in
  \mathbb{N}$, the class of matrix-valued functions
$$R(\cdot;\varepsilon):[0,\varepsilon_0]\rightarrow L_1 ^{m\times m}$$
parametrized by $ \varepsilon $
such that the solution of the Cauchy problem
$$
Z'(t;\varepsilon)= R( t;\varepsilon)Z(t;\varepsilon), \quad Z(a;\varepsilon) = I_m,
$$
satisfies the limit condition
$$
\lim\limits_{\varepsilon \rightarrow 0+} \|Z(\cdot;\varepsilon) - I_m\|_C =0,
$$
where $\|\cdot\|_C$ is the sup-norm.
\end{definition}

In paper \cite{MR2}, the following general result is established:

\begin{theorem}\label{1 limit G}
  Suppose that the vector boundary-value problem
\begin{equation}\label{bound probl 1}
    y'(t;\varepsilon)=A(t;\varepsilon)y(t;\varepsilon)+f(t;\varepsilon),\quad
t \in \mathcal{J}, \quad \varepsilon \in [0, \varepsilon_0],
\end{equation}
\begin{equation}\label{bound probl 2}
    U_{\varepsilon}y(\cdot;\varepsilon)= 0,
\end{equation}
where the matrix-valued functions $A(\cdot,\varepsilon) \in
L_{1}^{m\times m}$, the vector-valued functions
$f(\cdot,\varepsilon) \in L_{1}^{m}$, and the linear continuous
operators
$$U_{\varepsilon}:C(\overline{\mathcal{J}};\mathbb{C}^{m}) \rightarrow\mathbb{C}^{m}, \quad m \in \mathbb{N},$$
satisfy the following conditions.
\begin{align*}
  1) \quad &\text{The homogeneous limit boundary-value
    problem~}(\ref{bound probl 1}), (\ref{bound probl 2})
  \text{~with~} \varepsilon = 0  \text{~and~} f(\cdot;0) \equiv 0 \\
  &\text{~ has only a trivial solution;}\\
  2)\quad &A(\cdot;\varepsilon)-A(\cdot;0)\in \mathcal{M}^m;\\
  3)\quad &\|U_{\varepsilon} - U_{0}\|\rightarrow 0,\quad
  \varepsilon\rightarrow 0+.
\end{align*}
Then, for a small enough $\varepsilon$, there exist Green matrices
$G(t, s; \varepsilon)$ for problems (\ref{bound probl 1}),
(\ref{bound probl 2}) and, on the square $\mathcal{J}\times
\mathcal{J}$,
\begin{equation}\label{G}
 \|G(\cdot,\cdot;\varepsilon)-G(\cdot,\cdot;0)\|_\infty \rightarrow 0,\quad \varepsilon\rightarrow 0+,
\end{equation}
where $\|\cdot\|_\infty$ is the norm in the space $L_\infty$.
\end{theorem}

\begin{remark}
  Condition 3) in Theorem \ref{1 limit G} cannot be replaced with
  the weaker condition on the operator $U_{\varepsilon}$ to strongly
  converge, $U_{\varepsilon} \stackrel{s}\rightarrow U_{0}$
  \cite{MR2}.  However, one can easily see that, for the two-point
  boundary operators
$$U_\varepsilon y := B_1(\varepsilon) y(a) + B_2(\varepsilon) y(b),
\quad B_k(\varepsilon) \in \mathbb{C}^{m\times m},\quad k = 1,2,$$
both the strong convergence and the norm convergence conditions are
equivalent to
$${\|B_k(\varepsilon) - B_k(0)\| \rightarrow
0,}\quad \varepsilon\rightarrow 0+, \quad k = 1,2.$$
\end{remark}

There are different sufficient conditions for the matrix-valued
function $R(\cdot;\varepsilon)$ to belong to $\mathcal{M}^{m}$. In
particular, the results of \cite{Tamar} give that conditions (1),
(2), (3) of Theorem~\ref{res conv part} imply
$$
A(\cdot;\varepsilon)-A(\cdot;0)\in \mathcal{M}^2,
$$
where the matrix-valued function $A(\cdot;\varepsilon)$ is given by
the formula
\begin{equation}\label{A matrix}
A(\cdot;\varepsilon):=\left ( \begin{array}{cc}
Q_\varepsilon/p_\varepsilon& 1/p_\varepsilon\\
-Q^2_\varepsilon/p_\varepsilon& -Q_\varepsilon/p_\varepsilon
\end{array}\right) \in L_1^{2\times 2}.
\end{equation}
Before proving Theorem \ref{res conv part}, we will need the
following two lemmas require to reduce Theorem~\ref{res conv part}
to Theorem \ref{1 limit G}.

\begin{lemma}\label{lemm1}
  The function $y(t)$ is a solution of the boundary-value problem
\begin{equation}\label{D^2}
 l_\varepsilon[y](t)= f(t;\varepsilon) \in L_2 ,\quad\varepsilon\in [0,\varepsilon_0],
\end{equation}
\begin{equation}\label{alpha+beta}
  \alpha(\varepsilon)\mathcal{Y}_\varepsilon(a)+
  \beta(\varepsilon)\mathcal{Y}_\varepsilon(b)=0,
\end{equation}
if and only if the vector-valued function $w(t) = (y(t),
D^{[1]}_\varepsilon y(t))$ is a solution of the boundary-value
problem
\begin{equation}\label{diff eq}
w'(t)=A(t;\varepsilon)w(t) + \varphi(t;\varepsilon),
\end{equation}
\begin{equation}\label{diff alpha+beta}
\alpha(\varepsilon)w(a)+
  \beta(\varepsilon)w(b)=0,
  \end{equation}
  where the matrix-valued function $A(\cdot;\varepsilon)$ is given
  by (\ref{A matrix}) and $\varphi(\cdot;\varepsilon) := (0,
  -f(\cdot;\varepsilon))$.
\end{lemma}

\begin{proof}[Proof of Lemma \ref{lemm1}] Consider the system of equations
$$\left\{
\begin{aligned}
  ( D^{[0]}_\varepsilon y(t))' & = \frac{Q_\varepsilon(t)}{p_\varepsilon(t)}D^{[0]}_\varepsilon y(t) +
    \frac{1}{p_\varepsilon(t)}D^{[1]}_\varepsilon y(t),\\
   ( D^{[1]}_\varepsilon y(t))' & =
   - \frac{Q_\varepsilon^2(t)}{p_\varepsilon(t)}D^{[0]}_\varepsilon y(t) -
    \frac{Q_\varepsilon(t)}{p_\varepsilon(t)}D^{[1]}_\varepsilon y(t) - f(t; \varepsilon). \\
\end{aligned}
\right. $$

Let $y(\cdot)$ be a solution of (\ref{D^2}), then the definition of
a quasi-derivative implies that $y(\cdot)$ is a solution of this
system.  On the other hand, denoting $w(t) = (D^{[0]}_\varepsilon
y(t), D^{[1]}_\varepsilon y(t))$ and $\varphi(t;\varepsilon) = (0,
-f(t;\varepsilon))$, we rewrite this system in the form of equation
(\ref{diff eq}).

Taking into account that $\mathcal{Y}_\varepsilon(a) = w(a)$,
$\mathcal{Y}_\varepsilon(b) = w(b)$, one can see that the boundary
conditions (\ref{alpha+beta}) are equivalent to the boundary
conditions (\ref{diff alpha+beta}).
\end{proof}

Due to Lemma \ref{lemm1}, that statement that
\begin{itemize}
\item [$(\mathcal{U})$] the homogeneous boundary-value problem
  $l_0[y](t)=0, \quad \alpha(0)\mathcal{Y}_0(a)+
  \beta(0)\mathcal{Y}_0(b) = 0$, has only a trivial
  solution
\end{itemize}
implies that the homogeneous boundary-value problem
$$w'(t)=A(t;0)w(t), \quad \alpha(0)w(a) + \beta(0)w(b)=0$$
has only a trivial solution.

\begin{lemma}\label{Gamma exist}
Let a Green
matrix
 $$
 G(t,s,\varepsilon)=(g_{ij}(t,s))_{i,j=1}^2\in L_\infty^{2\times
   2}
 $$
 exist for the problem (\ref{diff eq}), (\ref{diff alpha+beta}) for
 small enough $\varepsilon$. Then there exists a Green function
 $\Gamma(t,s;\varepsilon)$ for the semi-homogeneous boundary-value
 problem (\ref{D^2}), (\ref{alpha+beta}) and
$$ \Gamma(t,s;\varepsilon) = -g_{12}(t,s;\varepsilon)
\quad\mbox{a.e.}$$
\end{lemma}

\begin{proof}[Proof of Lemma \ref{Gamma exist}]
  According to the definition of a Green matrix, a unique solution
  of problem (\ref{diff eq}), (\ref{diff alpha+beta}) can be written
  in the form
$$w_\varepsilon(t)=\int\limits_a ^b
G(t,s;\varepsilon)\varphi(s;\varepsilon) ds, \quad t\in \mathcal{J}.$$

Due to Lemma \ref{lemm1}, the latter equality can be rewritten in the
form
$$\left\{
\begin{aligned}
    D^{[0]}_\varepsilon y_\varepsilon(t) & = \int\limits_a^b g_{12}(t,s;\varepsilon)(-f(s;\varepsilon))ds, \\
    D^{[1]}_\varepsilon y_\varepsilon(t) & = \int\limits_a^b g_{22}(t,s;\varepsilon)(-f(s;\varepsilon))ds, \\
\end{aligned}
\right.$$ where $y_\varepsilon(\cdot)$ is a unique solution of the
problem (\ref{D^2}), (\ref{alpha+beta}).  This implies the statement
of Lemma \ref{Gamma exist}.
\end{proof}

\begin{proof}[Proof of Theorem \ref{res conv part}]
Note that, due to the equality
$$(Q_\varepsilon + \mu)^2/p_\varepsilon - (Q_0 + \mu)^2/p_0 = (Q_\varepsilon^2/p_\varepsilon - Q_0^2/p_0)
+ 2\mu(Q_\varepsilon/p_\varepsilon - Q_0/p_0) +
\mu^2(1/p_\varepsilon - 1/p_0),$$ where $\mu \in \mathbb{C}$,
conditions (1)--(3) of Theorem \ref{res conv part} imply that we can
assume without loss of generality that $0 \in \rho(L_0)$.

We need to show that $\sup\limits_{\|f\|_2 =
1} \|L_\varepsilon^{-1}f - L_0^{-1}f\| \rightarrow 0$, $\varepsilon
\rightarrow 0+$.

The equation $L_\varepsilon^{-1}f = y_\varepsilon$ is equivalent to
$L_\varepsilon y_\varepsilon = f$, i. e., $y_\varepsilon$ is a
solution of problem (\ref{D^2}), (\ref{alpha+beta}). Also the
statement ($\mathcal{U}$) is verified due to $0 \in \rho(L_0)$.
From the conditions 1)--3) of Theorem \ref{res conv part}, it
follows that $A(\cdot;\varepsilon)-A(\cdot;0)\in \mathcal{M}^2$,
where $A(\cdot;\varepsilon)$ is given by formula (\ref{A matrix}).
Thus statement of Theorem \ref{res conv part} implies that the
problem (\ref{diff eq}), (\ref{diff alpha+beta}) satisfies
conditions of Theorem \ref{1 limit G}. This means that Green
matrices $G(t,s;\varepsilon)$ of the problems (\ref{diff eq}),
(\ref{diff alpha+beta}) exist and the limit relation (\ref{G}) is
satisfied.  Taking into account Lemma \ref{Gamma exist} this yields
limit equality
 $$
 \|\Gamma(\cdot,\cdot;\varepsilon)-\Gamma(\cdot,\cdot;0)\|_\infty
 \rightarrow 0,\quad \varepsilon\rightarrow 0+.$$

Then
\begin{align*}
\|L_\varepsilon^{-1}- L_0^{-1}\| & = \sup\limits_{\|f\|_2=1}
\|\int_a^b\left[\Gamma(t,s;\varepsilon) -
\Gamma(t,s;0)\right]f(s)\,ds\|_2
\\
& \leq(b - a)^{1/2}\sup\limits_{\|f\|_2=1} \|\int_a^b
\left|\Gamma(t,s;\varepsilon) - \Gamma(t,s;0)\right|
\left|f(s)\right|ds\|_C
\\ &\leq (b -
a)\|\Gamma(\cdot,\cdot;\varepsilon) -
\Gamma(\cdot,\cdot;0)\|_\infty \rightarrow 0, \quad \varepsilon
\rightarrow 0+,
\end{align*}
which proves Theorem \ref{res conv part}.
\end{proof}

For the case $p_\varepsilon(t) \equiv 1$, a statement stronger than Theorem \ref{res conv part} was proved in
\cite{GM res conv}.

\section{Self-adjoint boundary conditions}
In what follows we will require the functions $p$, $Q$ and,
consequently, the distribution $q = Q'$ to be real-valued. In this
case, the expression $l[y]$ is formally self-adjoint \cite{EM} and,
according to Theorem \ref{L adjoint}, the minimal operator
$L_{\operatorname{min}}$ is symmetric. So one may pose a problem of
describing (in terms of homogeneous boundary conditions) all
extensions of the operator $L_{\operatorname{min}}$ that are
self-adjoint in the space $L_2$. To give an answer to this question,
we will apply the concept of the boundary triplet.

Let us recall following definition.
\begin{definition}
  Let $L$ be a closed densely defined symmetric operator in a Hilbert space
  $\mathcal{H}$ with equal (finite or infinite) deficient indices.
  The triplet $\left( {H,\Gamma _1 ,\Gamma _2 }\right)$, where $H$
  is an auxiliary Hilbert space and $\Gamma_1$, $\Gamma_2$ are the
  linear mappings of $\operatorname{Dom}(L^*)$ onto $H,$ is called
  a \emph{boundary triplet} of the symmetric operator $L$, if
\begin{enumerate}
  \item  for any $ f,g \in \operatorname{Dom} \left( {L^*} \right)$,
  $$
\left( {L^ *  f,g} \right)_\mathcal{H} - \left( {f,L^ *  g}
\right)_\mathcal{H} = \left( {\Gamma _1 f,\Gamma _2 g} \right)_H  -
\left( {\Gamma _2 f,\Gamma _1 g} \right)_H,
$$
\item for any $ f_1, f_2 \in H$ there is a vector $ f\in
  \operatorname{Dom} \left( {L^*} \right)$ such that $\Gamma _1 f =
  f_1$, $ \Gamma _2 f = f_2$.
\end{enumerate}
\end{definition}

The definition of a boundary triplet implies that $ f \in
\operatorname{Dom} \left( {L} \right)$ if and only if $\Gamma_1f =
\Gamma_2f = 0$. A boundary triplet exists for any symmetric operator
with equal non-zero deficient indices (see \cite{Gorbachuk} and
references therein). It is not unique.

The following result is crucial for the rest of the paper.

\begin{baslem*}
  \label{PGZth} Triplet $(\mathbb{C}^{2}, \Gamma_1, \Gamma_2)$,
  where $\Gamma_1, \Gamma_2$ are the linear mappings
\begin{equation} \label{PGZ}
\Gamma_1y := \left( D^{[1]}y(a), -D^{[1]}y(b)\right), \, \Gamma_2y := \left( y(a), y(b)\right),
\end{equation}
from $\operatorname{Dom}(L_{\operatorname{max}})$ onto
$\mathbb{C}^{2}$ is a boundary triplet for the operator
$L_{\operatorname{min}}$.
\end{baslem*}

For convenience, we introduce the following notation.
\begin{definition}
  Denote by $L_K$ the restriction of the operator
  $L_{\operatorname{max}}$ onto the set of functions ${y(t) \in
    \operatorname{Dom}(L_{\operatorname{max}})}$ satisfying the
  homogeneous boundary condition in the canonical form
\begin{equation} \label{rozsh}
 \left( {K - I} \right)\Gamma _1 y +
i\left( {K + I} \right)\Gamma _2 y = 0,
\end{equation}
where $K$ is any bounded operator on the space $\mathbb{C}^{2}$.
\end{definition}

Basic Lemma together with results of \cite[Ch.~3]{Gorbachuk} gives
the following description of all self-adjoint extensions of
$L_{\operatorname{min}}$.

\begin{theorem}\label{adj ext}
  Every $L_K$, with $K$ being a unitary operator on the space
  $\mathbb{C}^{2}$, is a self-adjoint extension of the operator
  $L_{\operatorname{min}}$.  Conversely, for any self-adjoint
  extension $\widetilde{L}$ of the operator $L_{\operatorname{min}}$
  there is a unitary operator $K$ such that $\widetilde{L} = L_K$.
  This correspondence between unitary operators $\{K\}$ and
  self-adjoint extensions $\{\widetilde{L}\}$ is bijective.
\end{theorem}

We start a proof of the Basic Lemma with the following two lemmas
that are special cases of the corresponding results for general
quasi-differential expressions (see \cite{EM}).

\begin{lemma}\label{Lagrange}
Suppose $y, z \in \operatorname{Dom}(L_{\operatorname{max}})$.
Then
 $$ \int\limits_a^b \left( D^{[2]}y\cdot\overline z  -
   y\cdot\overline{D^{[2]}z} \right) dt = \left.\left( -
     D^{[0]}y\cdot\overline{D^{[1]}z} +
     D^{[1]}y\cdot\overline{D^{[0]}z}\right)\right|^b_a. $$
\end{lemma}

\begin{lemma}\label{surject}
  Suppose that $ \{\alpha _0 ,\alpha _1 \}$, $\{ \beta _0,\beta _1
  \}$ are arbitrary sets of complex numbers. Then there is a
  function ${y \in \operatorname{Dom}(L_{\operatorname{max}})}$ such
  that
  $$ D^{[k]}y(a) = \alpha _k , \quad D^{[k]}y(b) = \beta _k,
  \quad k = 0,1.$$
\end{lemma}

\begin{proof}[Proof of the Basic Lemma]
  To prove the Basic Lemma, we need to prove that the triplet
  $(\mathbb{C}^{2}, \Gamma_1, \Gamma_2)$ satisfies conditions $1)$
  and $2)$ in the definition of the boundary triplet for the
  operator $L_{\operatorname{min}}$. According to Theorem \ref{L
    adjoint}, $L^*_{\operatorname{min}} = L_{\operatorname{max}}$.
  Due to Lemma \ref{Lagrange},
 $$\left( {L_{\operatorname{max}}y,z}
  \right) - \left( y,L_{\operatorname{max}}z\right) =
  \left.\left(D^{[0]}y\cdot\overline{D^{[1]}z} -
      D^{[1]}y\cdot\overline{D^{[0]}z}\right)\right|^b_a.$$
  But
\begin{align*}
\left( {\Gamma _1 y,\Gamma _2 z} \right) & =
D^{[1]}y(a)\cdot\overline{D^{[0]}z(a)} -
D^{[1]}y(b)\cdot\overline{D^{[0]}z(b)},\\ \left({\Gamma _2
y,\Gamma _1 z} \right) & = D^{[0]}y(a)\cdot\overline{D^{[1]}z(a)}
- D^{[0]}y(b)\cdot\overline{D^{[1]}z(b)}.
\end{align*}
This means that condition $1)$ is fulfilled. Condition $2)$ is true due to Lemma \ref{surject}.
\end{proof}

\begin{proof}[Proof of Theorem \ref{adj ext}]
  The claim in Theorem \ref{adj ext} follows from the Basic Lemma
  and Theorem~1.6 Ch.~3 \cite{Gorbachuk} for the boundary triplet of
  an abstract symmetric operator.
\end{proof}

\begin{remark}\label{homeo}
  Theorem \ref{res conv part}, together with Theorem \ref{adj ext},
  implies that the mapping $K \to L_K$ is not only bijective but
  also continuous. More accurately, if  unitary operators $K_n$
  converge to an operator $K$, then
  $$ \left\|\left(L_K -
      \lambda\right)^{-1} - \left(L_{K_n} -
      \lambda\right)^{-1}\right\| \rightarrow 0, \quad n \rightarrow
  \infty, \quad \operatorname{Im} \lambda \neq 0. $$ The converse is
  also true, because the set of unitary operators in the space
  $\mathbb{C}^2$ is a compact set. This means that the mapping
  $$K
  \to \left(L_K - \lambda\right)^{-1}, \quad \operatorname{Im}
  \lambda \neq 0,$$ is a homeomorphism for any fixed $\lambda \in
  \mathbb{C}\setminus\mathbb{R}$.
\end{remark}

Now we pass to a description of separated self-adjoint boundary
conditions for expression (\ref{vyraz}).

Denote by $\mathbf{f_a}$ the germ of a continuous function $f$ at the point $a$.

\begin{definition}
  The boundary conditions that define the operator $L \subset
  L_{\operatorname{max}}$ are called \emph{separated} if for
  arbitrary functions $y \in \operatorname{Dom}(L)$ and $g, h \in
  \operatorname{Dom}(L_{\operatorname{max}})$,
$$g, h \in \operatorname{Dom}(L) \quad \text{if} \quad
\mathbf{g_a} = \mathbf{y_a},\quad  \mathbf{g_b} = 0,\quad
\mathbf{h_a} = 0,\quad  \mathbf{h_b} = \mathbf{y_b}.$$
\end{definition}

\begin{theorem}\label{divided adj}
Self-adjoint boundary conditions (\ref{rozsh}) are separated
if and only if
the matrix $K$ is of the form (\ref{separable cond}),
where $K_a, K_b \in \mathbb{C}$ and $|K_a| = |K_b| = 1$.
\end{theorem}

A proof of Theorem \ref{divided adj} is based on the following
lemma.

\begin{lemma}\label{separated lemma}
  Boundary conditions of the form (\ref{rozsh}), with $K$ being any
  $\mathbb{C}^{2\times 2}$-matrix are separated if and only if
\begin{equation} \label{separable cond}
K = \left(%
\begin{array}{cc}
  K_a & 0 \\
  0 & K_b \\
\end{array}%
\right),
 \end{equation}
where $K_a, K_b \in \mathbb{C}$.
\end{lemma}

\begin{proof}[Proof of Lemma \ref{separated lemma}]
It is evident that $\mathbf{y_c} = \mathbf{g_c}$ implies
\begin{equation}\label{separated def}
y(c) = g(c), \quad (D^{[1]}y)(c) = (D^{[1]}g)(c), \quad c \in [a, b].
\end{equation}
Let the matrix $K$ have the form (\ref{separable cond}) in boundary
condition (\ref{rozsh}).  Then conditions (\ref{rozsh}) can be
written in the form of a system,
\[ \left\lbrace
\begin{aligned}
(K_a - 1)D^{[1]}y(a) + i(K_a + 1)y(a)& = 0,
\\
 -(K_b - 1)D^{[1]}y(b)
+ i(K_b + 1)y(b) & = 0.
\end{aligned}
 \right. \]

It is evident that these boundary conditions are separated.

Inversely, suppose that the boundary conditions (\ref{rozsh}) are
separated. The matrix $K \in \mathbb{C}^{2 \times 2}$ can be written
in the form
$$K = \left(%
\begin{array}{cc}
  K_{11} & K_{12} \\
  K_{21} & K_{22} \\
\end{array}%
\right).$$
We need to prove  $K_{12} = K_{21} = 0$.

Let us rewrite the boundary conditions (\ref{rozsh}) in the form of
the system
\[ \left\lbrace
\begin{aligned}
( K_{11} - 1)D^{[1]}y(a) - K_{12}D^{[1]}y(b) + i( K_{11} + 1)y(a)
+ iK_{12}y(b) & = 0,\\ K_{21}D^{[1]}y(a) - (K_{22} - 1)D^{[1]}y(b)
+ iK_{21}y(a) + i(K_{22} + 1)y(b)&  = 0.
\end{aligned}
 \right. \]

 The fact that the boundary conditions are separated implies that a
 function $g$ such that $\mathbf{g_a} = \mathbf{y_a}, \mathbf{g_b} =
 0$ also satisfies this system.  Due to equalities (\ref{separated
   def}) this gives
\[ \left\lbrace
\begin{aligned}
&K_{11} \left[D^{[1]}y(a) + iy(a)\right]= D^{[1]}y(a) - iy(a),\\
&K_{21}\left[D^{[1]}y(a) + iy(a)\right] = 0
\end{aligned}
 \right. \]
for any $y \in \operatorname{Dom}(L_K)$.

This means that either $K_{21} = 0$ or $D^{[1]}y(a) + iy(a) = 0$
for any $y \in \operatorname{Dom}(L_K)$. Suppose $K_{21} \neq 0$.

Let us return to the boundary conditions (\ref{rozsh}). For any $F =
(F_1, F_2) \in \mathbb{C}^2$, consider the vectors $-i \left(K +
  I\right)F$ and $\left(K - I\right)F$. Due to the Basic Lemma and
the definition of the boundary triplet, there exists a function $y_F
\in \operatorname{Dom}(L_{\operatorname{max}})$ such that
\begin{equation}\label{rozsh parametric}
\left\lbrace
\begin{aligned}
-i \left(K + I\right)F &  = \Gamma _1 y_F,\\ \left(K - I\right)F &
= \Gamma _2 y_F.
\end{aligned}\right.
\end{equation}
A simple calculation shows that $y_F$ satisfies the boundary
conditions (\ref{rozsh}) and, therefore, ${y_F \in}$
${\operatorname{Dom}(L_K)}$.  We can rewrite (\ref{rozsh
  parametric}) in the form of the system
\[
\left\lbrace
\begin{aligned}
- i(K_{11} + 1)F_1 - iK_{12}F_2 & = D^{[1]}y_F(a),
\\ -iK_{21}F_1 -
i(K_{22} + 1)F_2 & = -D^{[1]}y_F(b),
\\ (K_{11} - 1)F_1 + K_{12}F_2 & =
y_F(a),
\\ K_{21}F_1 + (K_{22} - 1)F_2 & = y_F(b).
\end{aligned}
\right.
\]

The first and the third equations of the system above imply that $0
= D^{[1]}y_F(a) + iy_F(a) = -2iF_1$ for any $F_1 \in \mathbb{C}$.
We arrived at a contradiction, therefore, $K_{21} = 0$.

Similarly one may prove $K_{12} = 0$.
\end{proof}

\begin{proof}[Proof of Theorem \ref{divided adj}]
  Due to Lemma \ref{separated lemma}, we only need to remark that
  a matrix of the form (\ref{separable cond}) is unitary if and
  only if $|K_a| = |K_b| = 1$.
\end{proof}

\section{Non-self-adjoint boundary conditions and generalized resolvents}
Recall the following definition.
\begin{definition}
  A densely defined linear operator $L$ on a complex Hilbert space
  $\mathcal{H}$ is called \emph{dissipative} if
 $$\operatorname{Im}
 \left( Lf, f \right)_\mathcal{H} \geq 0, \quad f \in
 \operatorname{Dom} (L) $$ and it is called \emph{maximal
   dissipative} if, besides this, $L$ has no nontrivial dissipative
 extensions on the space $\mathcal{H}$.
\end{definition}

For instance, every symmetric operator is dissipative and every
self-adjoint operator is a maximal dissipative one. Thus, if the
minimal operator $L_{\operatorname{min}}$ is symmetric, then one can
state the problem of describing its maximal dissipative extensions.
According to Phillips' Theorem \cite{Gorbachuk, Phil}, \emph{every
  maximal dissipative extension of a symmetric operator is a
  restriction of its adjoint operator.} Therefore, every maximal
dissipative extension of the operator $L_{\operatorname{min}}$ is a
restriction of operator $L_{\operatorname{max}}$.

Parametric bijective description of the class of maximal
dissipative extensions of the symmetric quasi-differential
operator $L_{\operatorname{min}}$ is given by the following
theorem.
\begin{theorem}
  \label{diss ext} Every $L_K$, with $K$ being a contracting
  operator on the space $\mathbb{C}^{2}$, is a maximal dissipative
  extension of the operator $L_{\operatorname{min}}$. Conversely,
  for any maximal dissipative extension $\widetilde{L}$ of the
  operator $L_{\operatorname{min}}$ there exists a contracting
  operator $K$ such that $\widetilde{L} = L_K$. This correspondence
  between contracting operators $\{K\}$ and the maximal dissipative
  extensions $\{\widetilde{L}\}$ is bijective.
\end{theorem}

\begin{proof}[Proof of Theorem \ref{diss ext}]
Theorem \ref{diss ext} is a direct consequence of Basic Lemma and
Theorem~1.6 Ch.~3 \cite{Gorbachuk} for the boundary triplet of an
abstract symmetric operator.
\end{proof}

\begin{remark}
The mapping
$$
K \rightarrow \left(L_K - \lambda\right)^{-1}, \quad
\operatorname{Im} \lambda < 0,
$$
for any fixed $\lambda$ is a homeomorphism (see Remark \ref{homeo}).
\end{remark}

\begin{theorem}\label{divided diss}
Dissipative boundary conditions (\ref{rozsh}) are separated
if and only if the matrix $K$ is of the form (\ref{separable cond}),
where $|K_a| \leq 1$, $|K_b| \leq 1$.
\end{theorem}

\begin{proof}[Proof of Theorem \ref{divided diss}]
As in the proof of Theorem \ref{divided adj}, due to Lemma \ref{separated lemma}, we only need to remark
that the matrix $K$ of the form (\ref{separable cond}) is a contracting operator on $\mathbb{C}^{2}$
if and only if
$|K_a| \leq 1$, $|K_b| \leq 1$.
\end{proof}

Recall the following definition.
\begin{definition}
  A \emph{generalized resolvent} of a closed symmetric operator $L$
  is the operator-valued function $R_\lambda$ of the complex parameter
  $\lambda \in \mathbb{C} \backslash \mathbb{R}$ which can be
  represented in the form
$$ R_\lambda
f = P^+ \left( L^+ - \lambda I^+\right)^{- 1}f, \quad f \in
\mathcal{H},$$ where $L^+$ is a self-adjoint extension of the
operator $L$, generally, on the space $\mathcal{H}^+$ which is wider
than $\mathcal{H}$, $I^+$ is the identity operator on
$\mathcal{H}^+,$ and $P^+$ is the orthogonal projection operator from
$\mathcal{H}^+$ onto $\mathcal{H}.$
\end{definition}

The operator-valued function $R_\lambda \quad (\operatorname{Im}
\lambda \neq 0)$ is a generalized resolvent of a symmetric
operator $L$ if and only if
 $$\left( R_\lambda f, g
 \right)_\mathcal{H} = \int_{-\infty}^{+\infty}\frac{d\left(F_\mu f,
     g\right)}{\mu - \lambda}, \quad f, g \in \mathcal{H}, $$
where
 $F_\mu$ is the generalized spectral function of the operator $L$. In other
 words, the operator-valued function $F_\mu$ should have following properties
 \cite{Ahiezer}:

 $1^0.$ For $\mu_2 > \mu_1$, the difference $F_{\mu_2} - F_{\mu_1}$
 is a bounded non-negative operator,

$2^0.$ $F_{\mu +} = F_\mu$ for any real $\mu$,

$3^0.$ for any $x \in \mathcal{H}$,
$$ \lim\limits_{\mu \rightarrow
- \infty}^{}||F_\mu x ||_\mathcal{H} = 0, \quad \lim\limits_{\mu
\rightarrow + \infty}^{} ||{F_\mu x - x} ||_\mathcal{H} = 0.$$

A parametric inner description of all generalized resolvents of the
operator $L_{\operatorname{min}}$ is given by the following theorem.

\begin{theorem}\label{gen res}
  There is a one-to-one correspondence between the generalized
  resolvents of the operator $L_{\operatorname{min}}$ and the
  boundary-value problems
  $$ l[y] = \lambda y + h,$$ $$ \left(
    {K(\lambda) - I} \right)\Gamma _1 y + i\left( {K(\lambda) + I}
  \right)\Gamma _2 y = 0, $$ where $\lambda \in \mathbb{C}$,
  $\operatorname{Im}\lambda < 0$, $h(x) \in L_2$, and $K(\lambda)$
  is an operator-valued function into the space $\mathbb{C}^2$,
  regular in the lower half-plane, such that $||K(\lambda)|| \leq
  1$.  This correspondence is given by the identity $$R_\lambda h =
  y,\quad \operatorname{Im}\lambda < 0.$$
\end{theorem}

\begin{proof}[Proof of Theorem \ref{gen res}]
Due to Basic Lemma Theorem, \ref{gen res} is a consequence of Theorem 1 of the paper \cite{Brook}.
\end{proof}

For general quasi-differential operators of even and odd orders,
respectively, the assertions of Theorems \ref{adj ext}, \ref{diss
  ext} and \ref{gen res} are announced without proofs in \cite{GM
  even, GM odd}.

\end{document}